\documentclass[12pt]{article} 
\usepackage{amsfonts,amsmath,amsxtra}
\usepackage{latexsym}
\usepackage{amssymb}

\def\hybrid{\topmargin 0pt      \oddsidemargin 0pt
        \headheight 0pt \headsep 0pt
        \textwidth 16.5cm
        \textheight 23cm
        \marginparwidth 0.0in
        \parskip 5pt plus 1pt   \jot = 1.5ex}
\catcode`\@=11
\def\marginnote#1{}
\newcount\hour
\newcount\minute
\newtoks\amorpm
\hour=\time\divide\hour by60 \minute=\time{\multiply\hour by60
\global\advance\minute by-\hour}
\edef\standardtime{{\ifnum\hour<12 \global\amorpm={am}%
        \else\global\amorpm={pm}\advance\hour by-12 \fi
        \ifnum\hour=0 \hour=12 \fi
      \number\hour:\ifnum\minute<10 0\fi\number\minute\the\amorpm}}
\edef\militarytime{\number\hour:\ifnum\minute<10 0\fi\number\minute}

\def\draftlabel#1{{\@bsphack\if@filesw {\let\thepage\relax
   \xdef\@gtempa{\write\@auxout{\string
      \newlabel{#1}{{\@currentlabel}{\thepage}}}}}\@gtempa
   \if@nobreak \ifvmode\nobreak\fi\fi\fi\@esphack}
        \gdef\@eqnlabel{#1}}
\def\@eqnlabel{}
\def\@vacuum{}
\def\draftmarginnote#1{\marginpar{\raggedright\scriptsize\tt#1}}

\def\draft{\oddsidemargin -0.1truein
        \def\@oddfoot{\sl preliminary draft \hfil
        \rm\thepage\hfil\sl\today\quad\militarytime}
        \let\@evenfoot\@oddfoot \overfullrule 3pt
        \let\label=\draftlabel
        \let\marginnote=\draftmarginnote
\def\@eqnnum{{\rm (\theequation)}
\rlap{\kern\marginparsep\tt\@eqnlabel}%
\global\let\@eqnlabel\@vacuum}  }

\newfont{\Bbbb}{msbm7 scaled 1\@ptsize00}
\newcommand{\zs}{\raise-1pt\hbox{$\mbox{\Bbbb Z}$}}

scaled 1\@ptsize00

\font\sevenmsa=msam6 
\newfam\msafam
\textfont\msafam=\sevenmsa
\def\hexnumber@#1{\ifnum#1<10 \number#1\else
\ifnum#1=10 A\else\ifnum#1=11 B\else\ifnum#1=12 C\else \ifnum#1=13
D\else\ifnum#1=14 E\else\ifnum#1=15 F\fi\fi\fi\fi\fi\fi\fi}
\def\msa@{\hexnumber@\msafam}
\def\llcorner{\delimiter"4\msa@78\msa@78 }
\def\lrcorner{\delimiter"5\msa@79\msa@79 }
\mathchardef\blacktriangleright="3\msa@49
\mathchardef\blacktriangleleft="3\msa@4A \font\tenmsb=msbm10 scaled
1\@ptsize00
\newfam\msbfam
\textfont\msbfam=\tenmsb \scriptfont\msbfam=\tenmsb

\newdimen\Squaresize \Squaresize=14pt
\newdimen\Thickness \Thickness=0.5pt

\def\Square#1{\hbox{\vrule width \Thickness
   \vbox to \Squaresize{\hrule height \Thickness\vss
      \hbox to \Squaresize{\hss#1\hss}
   \vss\hrule height\Thickness}
\unskip\vrule width \Thickness} \kern-\Thickness}

\def\Vsquare#1{\vbox{\Square{$#1$}}\kern-\Thickness}

\def\numberbysection{\@addtoreset{equation}{section}
        \def\theequation{\thesection.\arabic{equation}}}
\numberbysection

\renewcommand{\theequation}{\thesection.\arabic{equation}}
\def\titlepage{\@restonecolfalse\if@twocolumn\@restonecoltrue\onecolumn
     \else \newpage \fi \thispagestyle{empty}\c@page\z@
        \def\thefootnote{\fnsymbol{footnote}} }

\def\endtitlepage{\if@restonecol\twocolumn \else  \fi
        \def\thefootnote{\arabic{footnote}}
        \setcounter{footnote}{0}}  
\relax

\hybrid
\makeatletter
\newdimen\normalarrayskip            
\newdimen\minarrayskip               
\normalarrayskip\baselineskip \minarrayskip\jot
\newif\ifold             \oldtrue            \def\new{\oldfalse}
\def\arraymode{\ifold\relax\else\displaystyle\fi}
\def\eqnumphantom{\phantom{(\theequation)}} 
\def\@arrayskip{\ifold\baselineskip\z@\lineskip\z@
     \else
     \baselineskip\minarrayskip\lineskip1\baselineskip\fi}

\def\@arrayclassz{\ifcase \@lastchclass \@acolampacol \or
\@ampacol \or \or \or \@addamp \or
   \@acolampacol \or \@firstampfalse \@acol \fi
\edef\@preamble{\@preamble
  \ifcase \@chnum
     \hfil$\relax\arraymode\@sharp$\hfil
     \or $\relax\arraymode\@sharp$\hfil
     \or \hfil$\relax\arraymode\@sharp$\fi}}

\def\@array[#1]#2{\setbox\@arstrutbox=\hbox{\vrule
     height\arraystretch \ht\strutbox
     depth\arraystretch \dp\strutbox
width\z@}\@mkpream{#2}\edef\@preamble{\halign \noexpand\@halignto
\bgroup \tabskip\z@ \@arstrut \@preamble \tabskip\z@ \cr}%
\let\@startpbox\@@startpbox \let\@endpbox\@@endpbox
  \if #1t\vtop \else \if#1b\vbox \else \vcenter \fi\fi
  \bgroup \let\par\relax
  \let\@sharp##\let\protect\relax
  \@arrayskip\@preamble}
\def\eqnarray{\stepcounter{equation}%
              \let\@currentlabel=\theequation
              \global\@eqnswtrue
              \global\@eqcnt\z@
              \tabskip\@centering              
              \let\\=\@eqncr
              $$%
            \halign to \displaywidth  \bgroup
             \eqnumphantom \@eqnsel
      \hskip\@centering                               
    $\displaystyle  \tabskip\z@ {##}$%
    &\global\@eqcnt\@ne \hskip 2\arraycolsep
         $ \displaystyle  \arraymode{##}$\hfil
    &\global\@eqcnt\tw@ \hskip 2\arraycolsep
         $\displaystyle\tabskip\z@{##}$\hfil
         \tabskip\@centering
    &{##}\tabskip\z@\cr}
\makeatother

\def\IC{\mathbb{C}}

\def\IP{\mathbb{P}}

\def\CA {\mathcal{A}}

\def\CE {\mathcal{E}}
\def\CF {\mathcal{F}}
\def\CG {\mathcal{G}}

\def\CM {\mathcal{M}}
\def\CN {\mathcal{N}}
\def\CO {\mathcal{O}}

\def\CT {\mathcal{T}}

\def\apr {\overline {\partial }}

\def\Tr{{\rm Tr}}

\newcommand\bqa{\begin{eqnarray}}
\newcommand\eqa{\end{eqnarray}}
\def\be{\begin{eqnarray}\new\begin{array}{cc}}
\def\ee{\end{array}\end{eqnarray}}

\def\beq{\begin{equation}}
\def\eeq{\end{equation}}
\def\bse{\begin{subequations}}                
\def\ese{\end{subequations}}
\def\bp{\begin{pmatrix}}
\def\ep{\end{pmatrix}}

\def\stack#1#2{\raise0.7pt\hbox{$\mathrel{\mathop{#2}\limits^{#1}}$}}
\def\tr{\triangleright}
\def\tl{\triangleleft}
\def\sem{\mathsurround=0pt \raise1pt
\hbox{$\scriptscriptstyle>\!\!$}\:\!\!\tl}
\def\mes{\mathsurround=0pt \tr\!\:\!\raise0.8pt
\hbox{$\scriptscriptstyle\!\!<$}\,}
\def\]{\mathsurround=0pt ]\raise-2pt\hbox{$_\ast$}}

\def\<{\langle}
\def\>{\rangle}

\def\CO{{\cal O}}

\def\we{\raise-1pt\hbox{$\,\stackrel{\wedge}{,}\,$}}
\def\tr{{\rm tr}\,}
\def\Tr{{\rm Tr}\,}

\newcounter{pac}[section]

\newcounter{pacc}[subsection]

0

0

\begin{document}

\begin{flushright}
TCD-MATH-16-03
\end{flushright}

\vspace{15 mm}
\centerline{\Large \bf On spectral cover equations}
\vspace{3 mm}
\centerline{\Large \bf  in}
\vspace{3 mm}
\centerline{\Large \bf Simpson integrable systems}
\vspace{8 mm}
\centerline{Anton A. Gerasimov and Samson L. Shatashvili}
\vspace{2 mm}

\vspace{5 mm}
\begin{abstract}
\noindent {\bf Abstract}. Following Simpson we consider the integrable
system structure on the moduli spaces of Higgs bundles  on 
a compact K\"{a}hler manifold $X$.
We propose a description of the corresponding spectral cover of $X$  
as the fiberwise projective dual to a 
hypersurface  in the projectivization $\IP(\CT_{X} \oplus \CO_X)$ of the tangent
bundle $\CT_X$ to $X$. The defining equation of the hypersurface 
dual to the Simpson spectral
cover is explicitly  constructed in terms of the Higgs fields.   
\end{abstract}

\vspace{1cm}


\section{Introduction}

The notion  of a spectral curve, or, more precisely,  
a family of spectral curves   plays an important role in  the theory of
integrable systems. Thus the phase space of an integrable system
is interpreted as a total space of a family of
the Lagrangian abelian varieties identified with the Jacobian
varieties of the corresponding spectral curves. This provides effective means to
explicitly solve the integrable system via geometry of divisors on  
spectral curves. The standard construction of the families of
spectral curves is via characteristic polynomials of the Lax operators
depending on the spectral parameter.

In original constructions  the spectral 
parameter is a coordinate function on an underlying  elliptic curve
or its degeneration, and the spectral curve is a finite cover 
of the underlying curve realized
as a subvariety of its  cotangent bundle. 
More general construction was proposed by Hitchin \cite{H}.  
The phase space of the integrable system is 
identified with the moduli space of  Higgs bundles, a certain degeneration  
of the moduli space of complex $G$-bundles over an arbitrary algebraic
curve $\Sigma$. The corresponding 
spectral curve is given by a finite cover of $\Sigma$ realised 
as a hypersurface in the cotangent bundle $T^*\Sigma$, 
defined by the characteristic polynomial of the Higgs field on $\Sigma$. 

One of  the important applications of the spectral curve construction 
is to the  method of the quantum separation of variables   
 \cite{Sk1}, \cite {Sk2}. In this approach the equation defining
spectral curve as a hypersurface in $T^*\Sigma$  transforms into the
defining equation for corresponding 
quantum integrable system eigenfunctions  in the
separated variables. 
The Hitchin  integrable systems, both classical and quantum, are much studied and
have found many applications in representation theory,  quantum  field
theory and  string theory. 

A higher-dimensional generalisation of the  integrable
systems was proposed by Simpson \cite{S1}.
The phase space of the corresponding  integrable system is 
the moduli space of Higgs bundles  on  a  K\"{a}helr
manifold $X$ and the  role of the spectral curve is played by a
spectral cover of $X$ realized as a subvariety of the total space of some
vector bundle over $X$. A particularly interesting choice of the
vector bundle is given by the cotangent bundle $T^*X$ to $X$. 
By the analogy with the Hitchin description of the
spectral curves, it is possible to provide a set of polynomial equations 
describing the spectral cover of $X$ as a
subvariety of $T^*X$.
However, in the higher-dimensional  case the
corresponding set of equations turns out to be  highly overdetermined.
This might potentially be  a problem in the constructive approach  
to finding explicit solutions of the  corresponding integrable  system.  

In this note we point out a curious property of the Simpson spectral
cover. It is projectively dual to a hypersurface 
defined by a single explicit equation. 
 More precisely - consider the projective compactification
$\IP(\CT^*_X\oplus \CO)$ of the holomorphic cotangent bundle $\CT^*_X$.
We propose a description 
of the Simpson spectral cover in $\IP(\CT^*_X\oplus \CO)$ via 
fiberwise projective dual to a  hypersurface in 
the projective compactification $\IP(\CT_X\oplus \CO)$ of the tangent 
bundle $\CT_X$.  This hypersurface is defined 
by a single fiberwise polynomial function on   
$\IP(\CT_X\oplus \CO)$ (see equation \eqref{twsp}).
 In the case of $X$ being an algebraic curve this equation 
essentially coincides with the Hitchin equation of the spectral
curve. The formulation in terms of the tangent bundle instead   of
cotangent bundle bears a resemblance  with the  Legendre transform 
between the Hamiltonian and the Lagrangian
formulations of classical mechanics. 
One might hope that the proposed description  
of the Simpson spectral cover via projective duality will be useful in 
various applications of the theory 
of  integrable systems associated with higher dimensional algebraic
varieties.  

The theory of the Simpson integrable systems is an exciting area of
research expanding the traditional horizons of the theory of
integrable systems. For instance, the Simpson integrable systems,
similar to the Hitchin systems, allow a non-commutative deformation akin to 
the Knizhnik-Zamolodchikov connection leading to quantum versions of
the spectral covers. The corresponding
higher-dimensional analog of the quantum spectral curve shall be an
important part of the higher-dimensional holomorphic/chiral  quantum field
theory, which is an interesting direction to pursue.  
Note in this regard that the phase spaces of the 
Simpson integrable systems appear naturally in the description of 
$D$-brane string backgrounds and thus via dualities appear in various
related problems.  
One might also expect applications of the Simpson
integrable systems to some kind of 
higher dimensional version of the geometric Langlands
correspondence.

{\em Acknowledgments}: We are grateful to L. Katzarkov, A. Rosly and L. Takhtajan
for enjoyable  discussions.

\section{Hitchin integrable systems} 

We start  by  briefly  recalling  the construction of the
integrable systems associated with the holomorphic vector bundles over
algebraic curves due to Hitchin \cite{H}. 

Let $E$ be a rank $n$ complex $C^\infty$ vector bundle over a Riemann surface 
 $\Sigma$. Let $\Omega_\Sigma^{p,q}$ be 
  complex  vector bundles of differential forms of the type $(p,q)$
and $\Omega^{p,q}_\Sigma(E)=\Omega_\Sigma^{p,q}\otimes E$. To define 
a holomorphic vector bundle $\CE$ associated 
with $E$ it is enough to pick a differential operator $\apr_{\CE}:
\Omega^{0,0}_\Sigma(E)\to \Omega_\Sigma^{0,1}(E)$ 
such that $\apr_\CE(fs)=(\apr f)s+f\apr_{\CE}s$, for $s\in
\Omega^{0,0}_\Sigma(E)$ and $f\in C^{\infty}(\Sigma)$. 
Then the   holomorphic sections of $\CE$ satisfy the equation
$\apr_{\CE}s=0$. The space of complex structures on the
$C^\infty$-bundle $E$ is an 
infinite-dimensional affine space modeled on $\Omega^{0,1}_\Sigma({\rm
  End}(\CE))$.  The group $\CG$ 
of automorphisms of the vector bundle $E$ acts on 
the operator $\apr_{\CE}$ via conjugation $\apr_{\CE}\to
g^{-1}\apr_{\CE}g$. The moduli space $\CN$ of the stable complex
structures on a complex $C^{\infty}$ vector bundle $E$  
is obtained by taking the quotient of the space $\CA^s$ of the stable complex
structures on $E$  by the action of the automorphism group
$\CG$. Recall that the holomorphic vector bundle $\CE$ is stable iff 
for every proper subbundle $\CF\subset \CE$ one has 
$\frac{{\rm deg}(\CE)}{{\rm rank}(\CE)}>\frac{{\rm deg}(\CF)}{{\rm
    rank}(\CF)}$. 
In general the quotient space $\CN=\CA^s/\CG$ is  a smooth non-compact
space, but for the rank $n$ and the first Chern class $c_1(\CE)$
being  mutually prime the space  $\CN$ is compact. 

The cotangent bundle $T^*\CN$ to the moduli space $\CN$ can be
succinctly described via Hamiltonian reduction of the space of the Higgs  
bundles  with respect to the action of the group of automorphisms $\CG$. 
For a rank $n$ complex vector bundle $E$ let us consider the 
cotangent space $T^*\CA^s$ naturally identified with the space of  
pairs $(\apr_\CE,\Phi)$ of complex structures on $E$ defined by $\apr_\CE$ 
and $\Phi\in \Omega^{1,0}_{\Sigma}\otimes {\rm
End}(\CE)$. This space has a canonical holomorphic symplectic structure of the cotangent
bundle. The Hamiltonian reduction over zero set of moment map  with respect to
the action of the automorphism group $\CG$  provides a realization of the cotangent
bundle $T^*\CN$. Points of $T^*\CN$ parametrize
the Higgs bundles,  i.e. the pairs $(\apr_\CE,\Phi)$ satisfying the equation 
\be
 \apr_\CE \Phi=0.
\ee
There exists a more relaxed stability 
condition of the Higgs bundles (only $\Phi$-stable subbundles are
considered) leading to the Hitchin moduli space $\CM$
of stable Higgs bundles  such that $T^*\CN \subset \CM$ \cite{H}.

The total space of the cotangent bundle $T^*\CN$ is naturally a
holomorphic symplectic manifold with the holomorphic symplectic 
structure  descending from the $\CG$-invariant 
holomorphic symplectic structure 
\be\label{SS0}
\Omega=\int_{\Sigma_0}\,\,\,\Tr\,\delta \bar{A}\wedge \delta \Phi,
\ee 
on $T^*\CA^s$.  Here $\bar{A}\in \Omega^{0,1}_\Sigma({\rm End}(\CE))$
are affine coordinates on the space of complex structures on $\CE$,
i.e. for two complex structures  we have
$\apr^{(2)}_{\CE}=\apr^{(1)}_\CE+\bar{A}^{(2)}-\bar{A}^{(1)}$ and thus
the two-form \eqref{SS0} is well-defined. The holomorphic
  symplectic structure allows an extension on the moduli space $\CM$ of
  stable Higgs bundles. 

According to  Hitchin \cite{H}, 
the symplectic space $\CM$ is an  algebraically completely
integrable system. The corresponding set of mutually commuting algebraic functions
is provided by the Hitchin map 
\be\label{Hpr} 
\CM\longrightarrow 
\oplus_{i=1}^nH^0(\Sigma,{\rm Sym}^i(\Omega^{1,0}_{\Sigma})),
\ee   
sending the pair $(\apr_\CE,\Phi)$ to the vector $(\Tr \Phi,\Tr
\Phi^2, \cdots ,\Tr \Phi^n)$. To construct the explicit basis of the algebraic
Hamiltonian functions one shall fix a bases 
$\{\mu_{i,j}\}$, $i=1,\ldots, n$, $j=1,\ldots ,n_i$, 
$n_i=\dim H^1(\Sigma,{\rm Sym}^{i-1}\CT)$   
in the space  $\oplus_{i=1}^n H^1(\Sigma,{\rm Sym}^{i-1}\CT)$ 
where $\CT$ is the holomorphic tangent bundle to $\Sigma$.  
 Then the ring of  integrable 
 Hamiltonians is generated by the elements  
\be\label{HamH}
H_{i,j}=\int_{\Sigma} \mu_{i,j} \Tr \Phi^i,\qquad i=1,\ldots, n,, \qquad
j=1,\ldots ,n_i. 
\ee
Mutual commutativity of the Hamiltonians \eqref{HamH}  
follows from their invariance under the action of the automorphism
group $\CG$ and obvious mutual commutativity  
considered as functions on the space $T^*\CA^s$. 

Generic fibers of the  Hitchin map \eqref{Hpr} are  
compact  Lagrangian abelian varieties identified with 
the Jacobi varieties of the algebraic curves 
parametrized by the target space of the map \eqref{HamH}. 
This family of curves, known as a family of spectral curves of the integrable system,
 allows  the following  construction. 
One constructs a coherent sheaf
 $\hat{\CE}$ on the total space of $T^*\Sigma$ of the cotangent bundle
 to $\Sigma$ so that the spectral curve is given by the support of $\hat{\CE}$.
To define a coherent sheaf on $T^*\Sigma$ is 
the same as to define a sheaf of $\pi_* \CO_{T^*\Sigma}$-module 
on $\Sigma$. Here the direct image
$\pi_*\CO_{T^*\Sigma}$ is the sheaf of the polynomial algebras in one
variable over $\CO_\Sigma$ with the stalk  over a point $z\in \Sigma$
given by the algebra of the polynomial functions over the fiber
$T^*_z\Sigma$. Using the Higgs field $\Phi$ we can
supply $\CE$  with the structure of
$\pi_*\CO_{T^*\Sigma}$-module. For this it is enough to allow  
the generator of the polynomial algebra to act as multiplication by the
matrix $\Phi_z(z)$ where we locally have $\Phi=\Phi_z(z)dz$.   
Now the spectral curve  associated with   
$(\apr_\CE,\Phi)$ is given by the support of 
the coherent sheaf $\hat{\CE}$
on $T^*\Sigma$. Let $U\subset \Sigma$ be an open subset over which  the Higgs field
can be diagonalised with non-coincident eigenvalues. Then 
the construction above realises an open part of the spectral curve
as a non-ramified $n$-sheet cover of $U$. The spectral curve then
provides the  $n$-sheet cover of $\Sigma$ but, in general, ramified over
some points. It is clear that the 
 spectral curve, constructed this way, can be conveniently defined as 
the space of solutions of the characteristic equation of the Higgs field 
\be\label{HSC}
\det (w-\Phi(z))=0,  
\ee 
where $w$ is a linear coordinate on the fibers of the line bundle
$T^*\Sigma$. The coefficients of the spectral cover equation
\be
\det (w-\Phi(z))=\sum_{m=0}^n (-1)^{m}w^{n-m} \Tr_{\wedge^m} \Phi,
\ee
provide another set of generators of the ring of integrable  Hamiltonians 
\be\label{HamH1}
\tilde{H}_{m,j}=\int_{\Sigma} \mu_{m,j} \Tr_{\wedge^m} \Phi. 
\ee
Here we denote $\Tr_{\wedge^m} \Phi$ the trace of $\Phi$ in the $m$-th
fundamental representation of $GL_n$ realised in the space of $m$-th 
exterior powers of the standard representation $\IC^n$. 

At the end of this short review of the Hitchin construction let us 
note that the explicit description \eqref{HSC} 
of the spectral curve via the characteristic 
equation is widely used in the theory  of integrable systems and its
applications.

\section{ Simpson integrable systems}

A  generalization of the Hitchin integrable systems for
higher-dimensional base was introduced by Simpson \cite{S1}
using the notion of the Higgs bundles  in arbitrary dimensions.
In the following we consider a particular case of the Simpson's construction 
based on Higgs pairs  associated with the cotangent bundle of the
underlying manifold. 

Let $X$ be a $d$-dimensional 
K\"{a}hler manifold $X$ with the K\"{a}helr form $\omega$.  
Let $E$ be a rank $n$ complex vector bundle over $X$. The 
integrable complex structure on $E$ can be 
defined via  differential operator $\apr_{\CE}:
\Omega^{0,0}_X(E)\to \Omega^{0,1}_X(E)$ 
such that $\apr_\CE(fs)=(\apr f)s+f\apr_{\CE}s$, for $s\in
\Omega^{0,0}_X(\CE)$ and $f\in C^{\infty}(X)$ and  
the integrability condition $\apr_\CE^2=0$ holds.  Let $\CE$ be the
corresponding holomorphic vector bundle over $X$. 
The moduli space of  stable integrable complex structures modulo
action of the automorphism group $\CG$  of the vector bundle is 
a non-singular and  in general a non-compact K\"{a}helr manifold.

The generalized Higgs pair 
$(\apr_{\CE},\Phi)$ is a pair of: 1. an operator 
$\apr_{\CE}$ defining integrable complex structure on  $E$, and 2.
a holomorphic section $\Phi$ of $\Omega^{1,0}_{X}\otimes {\rm End}(\CE)$
which is squared to zero.  Thus we have the following defining
equations  for the Higgs pair $(\apr_{\bar{\CE}},\Phi)$
\be
\apr_{\CE}^2=0,\qquad \apr_{\CE}\Phi=0,\qquad \Phi\wedge
\Phi=0.  
\ee
Note that the last equation can be written in components as the
commutativity conditions  
\be\label{CC}
[\Phi_i,\Phi_j]=0, \qquad \Phi=\sum_{i=1}^d\Phi_i(z) dz^i. 
\ee
In the following we consider the case of $\CE$ having trivial rational
Chern classes. In this case to there is  one-to-one correspondence between
isomorphism classes of the Higgs bundles defined above  and the rank
$n$ complex local systems on $X$ \cite{S1}. 

Let $\CM$ be the  moduli space of stable Higgs bundles 
considered modulo  action of the automorphism group $\CG$. 
The stability condition on the Higgs bundles is given by the standard
stability condition  for vector bundles with the exception that only 
$\Phi$-stable subbundles  are considered \cite{S1}. 
The space $\CM$ is a non-singular
 holomorphic symplectic manifold with the 
holomorphic symplectic structure obtained by the reduction from the holomorphic
symplectic structure  
\be\label{HigS}
\Omega=\int_{X} \omega^{d-1}\,\Tr \delta \bar{A}\wedge \delta \Phi,   
\ee
on the space of pairs $(\apr_{\CE},\Phi)$ consisting of
not necessary integrable complex structures $\apr_{\CE}$
and $\Phi \in \Omega^{1,0}_{X}\otimes {\rm End}(\CE)$. 
Here $\omega$ is the K\"{a}hler structure on $X$ and $\bar{A}$ is the
affine coordinate on the space of not necessary integrable complex
structures. The reason for the existence of the holomorphic symplectic
structure (and moreover the hyperk\"{a}helr structure) is the
non-linear version of the Hodge decomposition on the first non-abelian
complex cohomology of $X$ (see \cite{S2} and references therein).   
The precise construction is analogous  to the construction of the holomorphic
symplectic structure (and more generally hyperk\"{a}hler structure) on
the moduli space of rank $n$ complex  local systems on $X$ (see \cite{F}
and references in \cite{S2}); moreover, the holomorphic symplectic
structure on the moduli
space of Higgs bundles can be obtained from the holomorphic structure
on the moduli space of complex local systems  by degeneration
using the identification of two moduli space  \cite{S1},
\cite{S2}. For a different approach see also \cite{DM}.

The moduli space $\CM$ is an algebraically completely integrable
system. The  analog of the Hitchin map \eqref{Hpr}  
\be\label{Higpr}
\CM \longrightarrow 
\oplus_{i=1}^n H^0(X,{\rm Sym}^i(\Omega^{1,0}_{X})),
\ee
sends the Higgs field $\Phi$ to the array $(\Tr \Phi,\Tr
\Phi^2, \cdots ,\Tr \Phi^n)$ of symmetric $GL_n$-invariant polynomial
functions. Precisely, let us contract the matrix valued one form 
$\Phi$ with a vector field $v$ on $X$ and consider traces $\Tr
(\iota_v\Phi)^i$ of its power. This gives degree $i$ symmetric
function of $v$ and thus an element of $H^0(X,{\rm Sym}^i(\Omega^{1,0}_{X}))$. 
The explicit set of Hamiltonians can be described using basis 
$\{\mu_{i,j}\}$, $i=1,\ldots, n$, $j=1,\ldots ,n_i$ 
in each  vector space $H^d(X,{\rm Sym}^{i-1}\CT\otimes \Omega_X^{d,0})$ 
where $\CT$ is the tangent bundle to $X$ 
(this is similar to Hitchin case, see Section 2).  
The generators of the ring of integrable 
 Hamiltonians are then given by 
\be\label{HamS}
H_{i,j}=\int_{X} \mu_{i,j} \Tr \Phi^i.  
\ee
The fibers of the projection \eqref{Higpr}
are Lagrangian abelian varieties that are Lagrangian with respect to
the symplectic structure on 
$\CM$. Simpson in \cite{S1}  defined an analog of the Hitchin spectral curve given
by a covering space $\tilde{X}\to X$. In good cases these spectral
covers are non-singular, and 
 the abelian varieties
in the fibration \eqref{Higpr} are identified with moduli spaces of
line bundles on $\tilde{X}$. 

These spectral covers can be represented as subvarieties 
of the cotangent space $T^*X$. Precisely they 
can be realized via supports of  coherent sheaf $\hat{\CE}$ on the total space of
$T^*X$ associated with the holomorphic vector bundle $\CE$
corresponding to the holomorphic structure on $E$ defined by $\apr_{\CE}$. 
This construction straightforwardly 
generalises  the corresponding construction of the spectral curve 
for the  Hitchin integrable system
described in the previous Section. To provide a description of
$\hat{\CE}$  we shall construct 
a sheaf of $\pi_* \CO_{T^*X}$-modules over $X$.
Note that $\pi_* \CO_{T^*X}$ is  a sheaf of polynomial algebras of
$d$-variables over $\CO_X$.  Using the Higgs field $\Phi$ we supply  $\CE$ 
with the structure of
$\pi_*\CO_{T^*X}$-module by allowing  the   
 generators of the polynomial algebra act as multiplications by the
matrices $\Phi_k(z)$ where we locally have
$\Phi=\sum_{k=1}^d\Phi_k(z)dz^k$. 
Now the spectral cover $\tilde{X}$  associated with   
$(\apr_\CE,\Phi)$ is given by the support of 
the coherent sheaf $\hat{\CE}$
on $T^*X$. An open part of the spectral cover is realised 
as  a  non-ramified $n$-sheet cover of 
 an open subset of $X$  over which  at
least one component $\Phi_k$ of the the Higgs field
can be diagonalised with non-coincident eigenvalues.
Let us give an explicit local construction of the
spectral covering locally over $X$.  
Let $(z^1,\ldots, z^d)$ be local coordinates over a small open subset
$U$ of $X$ and
$(w_1,\ldots , w_d)$ be local linear coordinates  in the fiber of the
cotangent bundle $T^*U$. Suppose that over $U$ the 
Higgs field  can be chosen in the diagonal form   
\be
\Phi=\sum_{k=1}^d\sum_{a=1}^nE_{aa}\Phi^a_k\,dz^k,
\ee
where $E_{ab}$ are elementary $(n\times n)$-matrices
$(E_{ab})_{\alpha\beta}=\delta_{a\alpha}\delta_{b\beta}$. Then the support of
$\hat{\CE}$  is locally a collection of $n$ copies of $U_0$ defined by
union of $n$ embeddings of $U$ into $T^*U$  
\be
(z^1,\ldots, z^d)\longrightarrow (z^1,\ldots,z^{d},\Phi_1^a(z),
\ldots \Phi_d^a(z)), \qquad a=1,\ldots, n. 
\ee
This provides a structure of degree $n$ non-ramified covering of $U$.   

The above construction of the spectral cover can be easily
reformulated in terms of a set of defining characteristic equations 
(see e.g. \cite{KOP}) generalizing the spectral curve description via
characteristic equation \eqref{HSC} for the Hitchin integrable systems. 
Let $(z^1,\ldots,z^d,w_1,\ldots, w_d)$
be local holomorphic coordinates on $T^*X$ and let $v$ be
a section of the holomorphic tangent bundle $\CT$. Then by
construction the linear operator $(\iota_v\Phi-\iota_vw)$ has non-trivial
kernel when $(z,w)\in \tilde{X}$ and thus $\det
(\iota_v\Phi-\iota_vw)=0$ on $\tilde{X}$. This shall hold for any
vector field $v$. Thus points of the spectral cover shall satisfy the following
equation 
\be\label{SSC}
\det(\Phi-w)=0,
\ee
as an element of $n$-th symmetric power of the cotangent bundle. 
This is a direct analog of \eqref{HSC}. 

The set of equation is overdetermined as the spectral cover 
is a codimension $d$ subspace in $T^*X$ while the condition 
\eqref{SSC} provides $\dim(S^n\IC^d)=\binom {n+d-1} {n}$-equations on
$T^*X$.  This renders  
the description \eqref{SSC} of the Simpson spectral cover less
explicit then the analogous equation \eqref{HSC} 
of the Hitchin spectral cover. 
In the following we demonstrate that
there is a convenient description of the Simpson spectral cover by a
single equation but in another space.   
This  construction uses fiberwise projective duality.

\section{Spectral cover via projective duality}

Let us start by  recalling  the standard projective duality relation
between  linear subspaces in complex projective spaces. Let  
$\IP(V)$, $\dim V=d+1$ be a projective space and $\IP(V^*)$ be its
dual. We fix dual homogeneous coordinates 
$(\xi^0,\ldots, \xi^{d})$ and $(\eta_0,\ldots , \eta_{d})$ 
 in $V$ and $V^*$ correspondingly.  
A linear $k$-dimensional subspace  $W_k \subset \IP(V)$ 
can be represented by a $(k+1)$-dimensional subspace of $V$ and thus   
be  described by a set
of $d-k$ linear homogeneous equations 
\be
\sum_{j=0}^{d}A^i_{j}\xi^j=0, \qquad  i=1,\ldots , d-k\,.
\ee
For instance the hyperplane in $\IP(V)$  
is a space of solutions of a linear homogeneous 
  equation 
\be\label{HPL}
\sum_{j=0}^{d}\alpha_{j}\xi^j=0.  
\ee
Linear subspaces generically intersect according to their
dimensions, i.e. a linear $k$-dimensional subspace $W_k$ in $\IP(V)$ 
intersects a linear  $m$-dimensional subspace $W_m$ 
over a linear subspace of the dimension $k+m-d$ for $k+m\geq d$.

The projective duality interchanges dimension $k$ linear space in $\IP(V)$,
$V=\IC^{d+1}$,  and dimension $(d-k-1)$ linear subspaces in $\IP(V^*)$.
Let us parametrize points of  the 
$k$-dimensional space $W_k\subset \IP(V)$  by 
 $(k+1)$-tuple $(s^0,\ldots ,s^{k})\in \IC^{k+1}$ of variables as follows  
\be
\xi^i=\sum_{j=0}^kA^{i}_js^j, \qquad i=0,\ldots , d.
\ee
The corresponding dual  $(d-k-1)$-linear subspace in $\IP(V^*)$ 
is defined by the set of linear homogeneous equations    
\be
\sum_{i=0}^{d}\eta_i A^{i}_j=0, \qquad  j=0,\ldots ,k. 
\ee
In particular  to the  hyperplanes $H$ in $\IP^d(V)$,   defined by the linear
equation \eqref{HPL}, corresponds the  dual point in $\IP(V^*)$  given by the
line in $V^*$
\be\label{PRII}
L_H=\{(\lambda \alpha_0, \ldots ,\lambda \alpha_d)|\lambda\in
\IC^*\}. 
\ee
This line  can be also characterized as a solution of a system of linear
homogeneous equations 
\be\label{PRIII}
\sum_{j=0}^d\beta_i^j\eta_j=0,\qquad i=1,\ldots ,d,
\ee
and \eqref{HPL} and \eqref{PRIII} are related by the consistency
condition 
\be\label{PRcon}
\sum_{j=0}^d\beta_i^j\alpha_j=0,\qquad i=1,\ldots ,d. 
\ee
The projective duality respects incidence relations between various linear
subspaces.  
For instance  two points belong to the same hyperplane in 
$\IP(V)$  are  dual to 
two hyperplanes in  $\IP(V^*)$ intersecting  at the point dual to the
hyperplane in $\IP(V)$. 

This classical projective duality picture allows various
generalisations (see e.g. \cite{GKZ}). 
In particular, one can define projective duals for non-linear subvarieties  
in $\IP(V)$. Precisely, the projective dual to a subvariety $Z\in \IP(V)$ is a
subvariety $Z^*$ in $\IP(V^*)$ given by the union of duals to all tangent
planes to $X$. More generally given a non-linear bundle of 
projective spaces over $X$  e.g. the  projectivisation $\IP(\CE_X\oplus
\CO)$ of a vector bundle $\CE_X$ over the base $X$ one can consider 
the bundle of the dual projective spaces $\IP(\CE_X^*\oplus \CO)$. In
this case the classical projective duality holds fiberwise.  

The main point of this note is the following statement: 

{\it Let $\IP(\CT_{X}\oplus \CO_{X})$
be the projectivization of the tangent bundle $\CT_{X}$ to a compact
complex manifold $X$,  $\dim X=d$,  
$(z^1,\ldots ,z^{d})$ be  local coordinates on $X$ and 
$(\xi^0,\ldots ,\xi^d)$ be homogeneous coordinates in the
fibers of the projective bundle $\IP(\CT_{X}\oplus \CO_{X})$.
Let $(\apr_\CE,\Phi)$ be a Higgs
pair associated with a rank $n$ holomorphic vector bundle $\CE$ on
$X$ and let the Higgs field  be locally  written as  
$\Phi=\sum_{j=1}^d\Phi_j(z)dz^j$. 
Consider the hypersurface $\widehat{X}$  in $\IP(\CT_{X}\oplus \CO_{X})$
 defined by the equation 
\be\label{twsp}
\det(\sum_{j=1}^d\xi^j\Phi_j(z)+\xi^0)=0.  
\ee
Then the fiberwise  projective dual to this hypersurface is the 
projective compactification of the Simpson spectral cover $\tilde{X}$ 
associated  with the Higgs  pair $(\apr_\CE,\Phi)$.}

Assume that a  subset $U\subset X$ is such that the Higgs field can be
diagonalized 
\be\label{Diag}
\Phi=\sum_{a=1}^n\sum_{j=1}^d E_{aa} \Phi^a_j(z) dz^j,
\ee
where $E_{ab}$ are elementary $(n\times n)$-matrices
$(E_{ab})_{\alpha\beta}=\delta_{a\alpha}\delta_{b\beta}$. 
Due to \eqref{CC} this is equivalent to the condition
that at least one of the Higgs field component is diagonalizable. 
Then the corresponding hypersurface \eqref{twsp} in 
$\IP(\CT_{X}\oplus \CO_{X})$
is given by the equation 
\be
\prod_{a=1}^n(\sum_{j=1}^d\Phi^a_{j}(z)\xi^j+\xi^0)=0, 
\ee
and thus is a family of unions of  $n$ hyperplanes  
\be
\sum_{j=1}^d\Phi^a_{j}(z)\xi^j+\xi^0=0, \qquad a=1,\ldots ,n, 
\ee
in the projective fibers of $\IP(\CT_{X}\oplus \CO_{X})$ 
  parametrized by $U\subset X$.  
Applying the projective duality correspondence \eqref{HPL},
\eqref{PRII} between hyperplanes and
points we obtain the family of $n$ points $a=1,\ldots ,n$ 
\be
(\eta_0,\eta_1,\ldots, \eta_d)=(\lambda,\lambda \Phi^a_1(z),\ldots ,\lambda \Phi^a_d(z)),
 \qquad  \lambda\in \IC^*, 
\ee
in the fibers of the dual
projective family $\IP(\CT_X^*\oplus \CO_X)$. 
This indeed defines the  spectral cover of $U\subset X$ as a subset of 
$\IP(\CT^*_{X}\oplus \CO_{X})$.    
We expect that this description can be compactified by adding
components where the Higgs field is not diagonalizable. This then   
provides a concise  projective dual description \eqref{twsp}
of the Simpson spectral cover.

Let us stress that the proposed construction \eqref{twsp} 
is trivially connected with the Hitchin spectral cover equation
\eqref{HSC}.  Following \eqref{twsp} the  
spectral cover $\widetilde{\Sigma}$ of the algebraic curve $\Sigma$ 
associated with the Higgs pair $(\apr_\CE,\Phi)$ 
is defined by the fiberwise homogeneous equation    
in $\IP(\CT_{\Sigma}\oplus \CO_{\Sigma})$
\be
\det(\xi^1\Phi_z(z)+\xi^0)=0, \qquad \Phi=\Phi_z(z)dz,  
\ee
and let us consider the diagonalized Higgs field \eqref{Diag} over a subspace
$U\subset \Sigma$.  
The dimension of the projective fiber is equal to one and thus the
projective duality maps points in the fibers of the projective bundle 
$\IP(\CT_{\Sigma}\oplus \CO_{\Sigma})$ into the points of the
projective fibers of $\IP(\CT^*_{\Sigma}\oplus \CO_{\Sigma})$ according
to the incidence relation \eqref{PRcon}. The incidence relation
\eqref{PRIII} is
easily solved so that the equation for the dual curve in 
$\IP(T^*_{\Sigma_0}\oplus \CO_{\Sigma_0})$
 is given by 
\be
\det(\eta_1-\eta_0\Phi_z(z))=0, 
\ee
which is indeed the projective compactification of the  Hitchin curve 
\be
\det(w-\Phi_z(z))=0.  
\ee
The Hitchin description \eqref{HSC} is recovered in the chart $\eta_0=1$, $\eta_1=w$.

\vskip 1cm

\noindent {\small {\bf A.G.}: {\sl Institute for Theoretical and
Experimental Physics, 117259, Moscow,  Russia; \hspace{8 cm}\,}

\noindent{\small {\bf S.Sh.}: {\sl  
School of Mathematics, Trinity College Dublin, Ireland; \hspace{8 cm}\,
\hphantom{xxxx}   \hspace{3 mm} Hamilton Mathematics Institute, TCD, Dublin 2, Ireland;
 \hspace{8 cm}\,
\hphantom{xxxx}   \hspace{3 mm} 
Chair Israel Gelfand, Institut des Hautes Etudes Scientifiques, France;
 \hspace{8 cm}\,
\hphantom{xxxx}   \hspace{3 mm} 
Simons Center for Geometry and Physics, Stony Brook, USA;
\vspace{1 mm}
 \hspace{8 cm}\,
\hphantom{xxxx}   \hspace{3 mm} 
\vspace{1 mm}
On leave of absence from: \hspace{16 cm}\,
\hphantom{xxxx}   \hspace{3 mm} 
Euler International Mathematics Institute, St. Petersburg, Russia
 \hspace{8 cm}\
\hphantom{xxxx}   \hspace{3 mm}
\vspace{1 mm} 
 \hspace{12 mm} 
Institute for Information Transmission Problems, Lab. 5, Moscow, Russia;
}\\
\hphantom{xxxxxx} {\it E-mail address}: {\tt samson@maths.tcd.ie}}


\begin{thebibliography}{AB}
\frenchspacing \smallbreak


\bibitem[DM]{DM} D. Donagi, E.~ Markman, {\it  
Spectral covers, algebraically completely integrable, hamiltonian
systems, and moduli of bundles}, Springer L.N.M 1620 (1996), 1--119. 

\bibitem[F]{F} A.~Fujiki, {\it   Hyper-K\"{a}hler structure on the
    moduli space of flat bundles}, Prospects in complex geometry 
(Katata and Kyoto, 1989), Springer L.N.M. 1468 (1991), 1--83.

\bibitem[GKZ] {GKZ} I. M.~Gelfand, M. M.~Kapranov, and A.V.~Zelevinsky,  
{\it Discriminants, Resultants, and Multidimensional Determinants},
Birkh\"{a}user Boston, 1994, 523 pp.

\bibitem[H]{H} N.~Hitchin, {\it Stable bundles and integrable
    systems},  Duke Math. J. 54 (1987), no. 1, 91--114. 


\bibitem[KOP]{KOP} L.~Katzarkov, D.~Orlov,  T.~Pantev {\it 
Notes on Higgs bundles and D-branes}, see webpage of CIMPA-CIMAT-SWAGP Workshop
 {\it Moduli Spaces and Mathematical Physics}, 21.01.2013-02.02.2013;

\bibitem[S1] {S1} C.~ Simpson, {\it 
Moduli of representations of the fundamental group of a smooth
projective variety.}, I: Publ. Math. I.H.E.S. 79 (1994), 47--129; 
II: Publ. Math. I.H.E.S. 80 (1994), 5--79.
 
\bibitem[S2]{S2} C.~Simpson, {\it  Non-abelian Hodge theory}.
 Proceedings, ICM-90, Kyoto Springer, Tokyo (1991), 198--230.

\bibitem[Sk1]{Sk1} E.K.~Sklyanin, 
{\it Separation of Variables in the Gaudin model}, 
Zapiski Nauchnykh Seminarov LOMI 164 (1987), 151–169, (in Russian); 
J. Sov. Math. 47 (1989), 2473–2488, (English transl.).

\bibitem[Sk2]{Sk2} E.K.~ Sklyanin, {\it Separation of variables. New
    trends}  in Quantum field theory, integrable models and beyond 
(Kyoto, 1994). Progr. Theoret. Phys. Suppl. 118 (1995), 35–60.




\end{thebibliography}
\end{document}